# How definitive *is* the standard interpretation of Gödel's Incompleteness Theorem?

Bhupinder Singh Anand

Standard interpretations of Gödel's "undecidable" proposition, [(A$x$)$R$($x$)], argue that, although [~(A$x$)$R$($x$)] is PA-provable if [(A$x$)$R$($x$)] is PA-provable, we may not conclude from this that [~(A$x$)$R$($x$)] is PA-provable. We show that such interpretations are inconsistent with a standard Deduction Theorem of first order theories.

## Contents

## 1. Introduction

In his seminal 1931 paper [Go31a], Gödel meta-mathematically argues that his "undecidable" proposition, [(A$x$)$R$($x$)][1], is such that (cf. [An02b], §1.6(*iv*)):

---

[1] We use square brackets to differentiate between a formal expression [$F$] and its interpretation "$F$", where we follow Mendelson's definition of an interpretation M of a formal theory K, and of the interpretation of a formula of K under M ([Me64], p49, §2).

If $[(Ax)R(x)]$ is PA-provable, then $[\sim(Ax)R(x)]$ is PA-provable.

Now, a standard Deduction Theorem of an arbitrary first order theory states that ([Me64], p61, Corollary 2.6):

If $T$ is a set of well-formed formulas of an arbitrary first order theory K, and if $[A]$ is a closed well-formed formula of K, and if $(T, [A])|\text{-}_K [B]$, then $T|\text{-}_K ([A => B])$.

In an earlier paper ([An02b], Appendix 1), we implicitly assumed, without proof, that:

$(T, [A])|\text{-}_K [B]$ holds if, and only if, $T|\text{-}_K [B]$ holds when we assume $T|\text{-}_K [A]$.   (*)

In other words, we assumed that $[B]$ is a deduction from $(T, [A])$ in K if, and only if, whenever $[A]^2$ is a hypothetical deduction from $T$ in K, $[B]$ is a deduction from $T$ in K.

We then argued, that it should follow (essentially by the reasoning in §2.2 below), that:

$[(Ax)R(x) => \sim(Ax)R(x)]$ is PA-provable,

and, therefore, that:

$[\sim(Ax)R(x)]$ is PA-provable.

We then concluded that PA is omega-inconsistent. However, these conclusions are inconsistent with standard interpretations of Gödel's reasoning, which, firstly, assert both $[(Ax)R(x)]$ and $[\sim(Ax)R(x)]$ as PA-unprovable, and, secondly, assume that PA can be omega-consistent. Such interpretations, therefore, implicitly deny that the PA-provability of $[\sim(Ax)R(x)]$ can be inferred from the above meta-argument; ipso facto, they imply that (*) is false.

---

[2] For the purposes of this paper, we assume everywhere that $[A]$ is a closed well-formed formula of K.

In the following sections, we review the Deduction Theorems used in the earlier argument, and give a meta-mathematical proof of (*). It follows that the standard interpretations of Gödel's reasoning are inconsistent with a standard Deduction Theorem of an arbitrary first order theory ([Me64], p61, Corollary 2.6). We conclude that such interpretations cannot be accepted as definitive.

### 1.1 An overview

We first review, in Theorem 1, the proof of a standard Deduction Theorem, if $(T, [A])|\text{-}_K [B]$, then $T|\text{-}_K [A => B]$, where an explicit deduction of $[B]$ from $(T, [A])$ is known.

We then show, in Corollary 1.2, that Theorem 1 can be constructively extended to cases where $(T, [A])|\text{-}_K [B]$ is established meta-mathematically, and where an explicit deduction of $[B]$ from $(T, [A])$ is not known.

We finally prove (*) in Theorem 2.

## 2. A standard Deduction Theorem

The following is, essentially, Mendelson's proof of a standard Deduction Theorem ([Me64], p61, Proposition 2.4) of an arbitrary first order theory K:

**Theorem 1**: If $T$ is a set of well-formed formulas of an arbitrary first order theory K, and if $[A]$ is a closed well-formed formula of K, and if $(T, [A])|\text{-}_K [B]$, then $T|\text{-}_K [A => B]$.

**Proof**: Let $<[B_1], [B_2], ..., [B_n]>$ be a deduction of $[B]$ from $(T, [A])$ in K.

Then, by definition, $[B_n]$ is $[B]$ and, for each $i$, either $[B_i]$ is an axiom of K, or $[B_i]$ is in $T$, or $[B_i]$ is $[A]$, or $[B_i]$ is a direct consequence by some rules of inference of K of some of the preceding well-formed formulas in the sequence.

We now show, by induction, that $T \vdash_K [A \Rightarrow B_i]$ for each $i =< n$. As inductive hypothesis, we assume that the proposition is true for all deductions of length less than $n$.

(i) If $[B_i]$ is an axiom, or belongs to $T$, then $T \vdash_K [A \Rightarrow B_i]$, since $[B_i \Rightarrow (A \Rightarrow B_i)]$ is an axiom of K.

(ii) If $[B_i]$ is $[A]$, then $T \vdash_K [A \Rightarrow B_i]$, since $T \vdash_K [A \Rightarrow A]$.

(iii) If there exist $j, k$ less than $i$ such that $[B_k]$ is $[B_j \Rightarrow B_i]$, then, by the inductive hypothesis, $T \vdash_K [A \Rightarrow B_j]$, and $T \vdash_K [A \Rightarrow (B_j \Rightarrow B_i)]$. Hence, $T \vdash_K [A \Rightarrow B_i]$.

(iv) Finally, suppose there is some $j < i$ such that $[B_i]$ is $[(Ax)B_j]$, where $x$ is a variable in K. By hypothesis, $T \vdash_K [A \Rightarrow B_j]$. Since $x$ is not a free variable of $[A]$, we have that $[(Ax)(A \Rightarrow B_j) \Rightarrow (A \Rightarrow (Ax)B_j]$ is PA-provable. Since $T \vdash_K [A \Rightarrow B_j]$, it follows by Generalisation that $T \vdash_K [(Ax)(A \Rightarrow B_j)]$, and so $T \vdash_K [A \Rightarrow (Ax)B_j]$, i.e. $T \vdash_K [A \Rightarrow B_i]$.

This completes the induction, and Theorem 1 follows as the special case where $i = n$. ¶[3]

## 2.1 A number-theoretic corollary

Now, Gödel has defined ([Go31a], p22, Definition 45(6)) a primitive recursive number-theoretic relation $xB_{(K, T)}y$ that holds if, and only if, $x$ is the Gödel-number of a deduction from $T$ of the K-formula whose Gödel-number is $y$.

We thus have:

**Corollary 1.1**[4]: If the Gödel-number of the well-formed K-formula $[B]$ is $b$, and that of the well-formed K-formula $[A \Rightarrow B]$ is $c$, then Theorem 1 holds if, and only if[5]:

---

[3] We use the symbol "¶" as an end-of-proof marker.

$(Ex)xB_{(K, T, [A])}b => (Ez)zB_{(K, T)}c$

## 2.2 An extended Deduction Theorem

We next consider the proposition:

**Corollary 1.2**: If we assume Church's Thesis[6], then Theorem 1 holds even if the premise $(T, [A])|-_K [B]$ is established meta-mathematically, and a deduction <$[B_1], [B_2], ..., [B_n]$> of $[B]$ from $(T, [A])$ in K is not known explicitly.

**Proof**: Since Gödel's number-theoretic relation $xB_{(K, T)}y$ is primitive recursive, it follows that, if we assume Church's Thesis - which implies that a number-theoretic relation is decidable if, and only if, it is recursive - we can effectively determine some finite natural number $n$ for which the assertion $nB_{(K, T, [A])}b$ holds, where the Gödel-number of the well-formed K-formula $[B]$ is $b$.

Since $n$ would then, by definition, be the Gödel-number of a deduction <$[B_1], [B_2], ..., [B_n]$> of $[B]$ from $(T, [A])$ in K, we may thus constructively conclude, from the meta-mathematically determined assertion $(T, [A])|-_K [B]$, that some deduction <$[B_1], [B_2], ..., [B_n]$> of $[B]$ from $(T, [A])$ in K can, indeed, be effectively determined.

Theorem 1 follows. ¶

---

[4] We note that Corollary 1.1 and Corollary 2.2 may be essentially different number-theoretic assertions, which may not be obviously equivalent; the "obvious" assumption (*), thus, may need a proof.

[5] We note that this is a semantic meta-equivalence, based on the definition of the primitive recursive relation $xB_{(K, T)}y$.

[6] Church's Thesis: A number-theoretic function is effectively computable if, and only if, it is recursive ([Me64], p147, footnote). We appeal explicitly to Church's Thesis here to avoid implicitly assuming that every recursive relation is algorithmically decidable (cf. [An02c], §II(7) Corollary 14.3). In Anand ([An02g], §2.5(xii)) we show that, under a constructive interpretation of classical foundational concepts, Church's Thesis is a Theorem; such a premise would not, then, be needed.

## 3. An additional deduction theorem

We finally prove (*) as an additional deduction theorem, in an arbitrary first order theory K:

**Theorem 2**: If K is an arbitrary first order theory, and if [A] is a closed well-formed formula of K, then $(T, [A])|-_K [B]$ if, and only if, $T|-_K [B]$ holds when we assume $T|-_K [A]$.

**Proof**: Firstly, if there is a deduction $<[B_1], [B_2], ..., [B_n]>$ of [B] from (T, [A]) in K, and there is a deduction $<[A_1], [A_2], ..., [A_m]>$ of [A] from T in K, then $<[A_1], [A_2], ..., [A_m], [B_1], [B_2], ..., [B_n]>$ is a deduction of [B] from T in K. Hence we have: if $(T, [A])|-_K [B]$, then $T|-_K [B]$ holds when we assume $T|-_K [A]$.

Secondly, if there is a deduction $<[B_1], [B_2], ..., [B_n]>$ of [B] from T in K, then we have, trivially, that: if $T|-_K [B]$ holds when we assume $T|-_K [A]$, then $(T, [A])|-_K [B]$.

Lastly, we assume that there is no deduction $<[B_1], [B_2], ..., [B_n]>$ of [B] from T in K. If, now, $T|-_{K'} [B]$ holds when we assume $T|-_{K'} [A]$ in any consistent extension K' of K, then, if we assume that there is a sequence $<[A_1], [A_2], ..., [A_m]>$ of well-formed K'-formulas such that $[A_m]$ is [A] and, for each $m >= i >= 1$, either $[A_i]$ is an axiom of K', or $[A_i]$ is in T, or $[A_i]$ is a direct consequence by some rules of inference of K' of some of the preceding well-formed formulas in the sequence, then we can show, by induction on the deduction length $n$, that there is a sequence $<[B_1], [B_2], ..., [B_n]>$ of well-formed K-formulas such that $[B_1]$ is $[A]$[7], $[B_n]$ is [B] and, for each $i > 1$, either $[B_i]$ is an axiom of K, or $[B_i]$ is in $T$, or $[B_i]$ is a direct consequence by some rules of inference of K of some of the preceding well-formed formulas in the sequence.

---

[7] [A] is thus the hypothesis in the sequence; it is the only well-formed K-formula in the sequence that is not an axiom of K, not in $T$, and not a direct consequence of the axioms of K by any rules of inference of K.

Hence, if there is a deduction $<[A_1], [A_2], ..., [A_m]>$ of $[A]$ from $T$ in K', then $<[A_1], [A_2], ..., [A_m], [B_2], ..., [B_n]>$ is a deduction of $[B]$ from $T$ in K'. By definition, it follows that $<[B_2], ..., [B_n]>$ is a deduction of $[B]$ from $(T, [A])$ in K. We thus have: if $T|\text{-}_K [B]$ holds when we assume $T|\text{-}_K [A]$, then $(T, [A])|\text{-}_K [B]$. This completes the proof. ¶

In view of Corollary 1.2, we thus have:

**Corollary 2.1**: If we assume Church's Thesis, and if $[A]$ is a closed well-formed formula of K, then we may conclude $T|\text{-}_K ([A] => [B])$ if $T|\text{-}_K [B]$ holds when we assume $T|\text{-}_K [A]$.[8]

We note that, in the notation of Corollary 1.1, if the Gödel-number of the well-formed K-formula $[A]$ is $a$, then Corollary 2.1 holds if, and only if[9]:

**Corollary 2.2**: $((Ex)xB_{(K, T)}a => (Eu)uB_{(K, T)}b) => (Ez)zB_{(K, T)}c$.

---

[8] We note that there is a model-theoretic proof of Corollary 2.1. The case $T|\text{-}_K [B]$ is straightforward.

If $\sim T|\text{-}_K [B]$, then, as noted in Theorem 2, if $T|\text{-}_K [B]$ holds when we assume $T|\text{-}_K [A]$, then there is a sequence $<[B_1], [B_2], ..., [B_n]>$ of well-formed K-formulas such that $[B_1]$ is $[A]$, $[B_n]$ is $[B]$ and, for each $i > 1$, either $[B_i]$ is an axiom of K, or $[B_i]$ is in $T$, or $[B_i]$ is a direct consequence by some rules of inference of K of some of the preceding well-formed formulas in the sequence.

(Note: In the following, if $T$ is the set of well-formed K-formulas $\{[T_1], [T_2], ..., [T_l]\}$ then $(T \& [A])$ denotes the well-formed K-formula $[T_1 \& T_2 \& ..., T_l \& A]$, and, $(T \& A)$ denotes its interpretation in M: $T_1 \& T_2 \& ..., T_l \& A$.)

If, now, any well-formed formula in $(T, [A])$ is false under an interpretation M of K, then $(T \& A) => B$ is vacuously true in M.

If, however, all the well-formed formulas in $(T, [A])$ are true under interpretation in M, then the sequence $<[B_1], [B_2], ..., [B_n]>$ interprets as a deduction in M, since the interpretation preserves the axioms and rules of inference of K (cf. [Me64], p57). Thus $[B]$ is true in M, and so is $(T \& A) => B$.

In other words, we cannot have $(T, [A])$ true and $[B]$ false in M as this would imply that there is some consistent extension K' of K in which $T|\text{-}_{K'} [A]$, but not $T|\text{-}_{K'} [B]$, which is contrary to the hypothesis that, in any consistent K in which we assume $T|\text{-}_K [A]$, we also have $T|\text{-}_K [B]$.

Hence, $(T \& A) => B$ is true in all models of K. By a consequence of Gödel's Completeness Theorem for an arbitrary first order theory ([Me64], p68, Corollary 2.15($a$)), it follows that $|\text{-}_K (T \& [A]) => [B])$, and, ipso facto, that $T|\text{-}_K ([A] => [B])$.

[9] We note that this, too, is a semantic meta-equivalence, based on the definition of the primitive recursive relation $xB_{(K, T)}y$.

## 4. Conclusion

Since standard interpretations of Gödel's reasoning and conclusions do not admit Theorem 2 as a valid inference, such interpretations are inconsistent with the standard Deduction Theorem for an arbitrary first order theory [Me64], p61, Proposition 2.4); they cannot, therefore, be considered definitive.

*Acknowledgement: My thanks to Hitoshi Kitada, Department of Philosophy, University of Tokyo, for pointing out that there is a distinction between the standard interpretation, and the author's earlier interpretation, of the standard Deduction Theorem. My thanks to Professor Leo Harrington for pointing out that the meaning of "$T|\text{-}_K [A] \Rightarrow T|\text{-}_K [B]$" was ambiguous in my original argument.*

*Author's e-mail: anandb@vsnl.com*

(*Updated: Monday 23$^{rd}$ June 2003 2:35:32 PM IST by re@alixcomsi.com*)